\documentclass[12pt]{article}
\usepackage{amsxtra,amssymb,amsthm,amsmath,latexsym}

\theoremstyle{plain}

\def\R{{\mathbb R}}

\def\C{{\mathbb C}}

\def\oH{{\overset{\circ}{H}}}
\def\oH1{{\overset{\circ}{H}\kern-.02in{}^1}}

\def\bee{\begin{equation*}}
\def\eee{\end{equation*}}
\def\be{\begin{equation}}
\def\ee{\end{equation}}

\begin{document}

This paper is published:

Alexander G. Ramm

Theory of hyper-singular integrals and its application to the Navier-Stokes problem, 

Contrib. Math. 2, (2020), 47–-54

Open access Journal:
www.shahindp.com/locate/cm;   

DOI: 10.47443/cm.2020.0041

\title{Theory of hyper-singular integrals and its application to the Navier-Stokes problem
}

\author{Alexander G. Ramm\\
 Department  of Mathematics, Kansas State University, \\
 Manhattan, KS 66506, USA\\
ramm@math.ksu.edu\\
http://www.math.ksu.edu/\,$\sim$\,ramm}

\date{}
\maketitle\thispagestyle{empty}

\begin{abstract}

In this paper the convolution integrals  $\int_0^t(t-s)^{\lambda -1}b(s)ds$ with hyper-singular kernels are considered,  where
 $\lambda\le 0$ and $b$ is a smooth or 
 $b$ is in $L^1(\R_+)$. For such $\lambda$
these integrals diverge classically even for smooth $b$.
These convolution integrals are defined in this paper for $\lambda\le 0$,
$\lambda\neq 0,-1,-2,...$. 

Integral equations and inequalities are considered with the hyper-singular 
kernels
$(t-s)^{\lambda -1}_+$ for $\lambda\le 0$, where $t^\lambda_+:=0$ for $t<0$.

 In particular, one is interested in the value $\lambda=-\frac 14$ because it is important for the Navier-Stokes problem (NSP). 
 
Integral equations of the type $b(t)=b_0(t)+
\int_0^t(t-s)^{\lambda-1}b(s)ds$, $\lambda\le 0$,
are studied.
 
 The solution of these equations is investigated,
existence and uniqueness of the solution is proved for $\lambda=-\frac 1 4$. This special value of $\lambda$ is of basic importance for a study of the Navier-Stokes problem (NSP).

The above results are applied to the analysis of the NSP
in the space $\R^3$ without boundaries. It is proved that the NSP is contradictory in
the following sense: even if one assumes that the initial data $v(x,0)>0$ one proves that the solution $v(x,t)$
to the NSP  has the property $v(x,0)=0$, in general. This paradox shows that the NSP is not a correct description of the fluid mechanics problem and it proves that the NSP does not have a solution, in general.

\end{abstract}

\section{Introduction }\label{S:1}
{\bf 1.1.} In this paper the new definition of the convolution
$\int_0^t(t-s)^{\lambda-1}b(s)ds$ with hyper-singular functions is 
given. We compare this definition with the one, based on the distribution theory, \cite{GS}. The $b(s)$ is assumed locally integrable function on $\R_+:=[0,\infty)$. This assumption is satisfied
in the Navier-Stokes problem (NSP), see \cite{R691},
Chapter 5, where the integral equations of the type
$b(t)=b_0(t)+\int_0^t (t-s)^{\lambda-1}b(s)ds$  with $\lambda=-\frac 1 4$ are of interest. Classically these integral equations do not make sense because the integrals diverge if $\lambda\le 0$. In Sections 1--4 of this paper a new definition of such integrals is given, the solution
to the integral equation with hyper-singular kernel is
 investigated. These results are used  in Section 5 of this paper, where
the basic results concerning the NSP are obtained.

We analyze the NSP and prove that the NSP is physically not a correct description of the fluid mechanics problem and that the NSP does not have a solution, in general.
The words "in general" mean that if the initial velocity $v_0=0$ and the force $f=0$, then the NSP has the solution
$v(x,t)=0$ which exists for all $t\ge 0$. This meaning
is valid in Theorem 5.2, below.

 For future use we define $\Phi_\lambda(t):=
\frac {t^{\lambda-1}}{\Gamma(\lambda)}$ and the convolution
$\Phi_\lambda \star b:=\int_0^t\Phi_\lambda(t-s) b(s)ds$.
Here and below $t:=t_+$, that is, $t=0$ if $t<0$ and $t=t$ if $t\ge 0$.
 
{\bf 1.2.}  Let us give the standard definition of the singular integral
used in the distribution theory. Let 
\be\label{e100}
 J:=\int_0^\infty t^{\lambda-1} \phi(t)dt,
 \ee
where the test function $\phi(t)\in C^\infty_0(\R)$.

Integral \eqref{e100} diverges classically (that is, in the classical sense) if $\lambda\le 0$. It is defined in distribution theory
(for example, in \cite{GS}) as follows: 
\be\label{e101}
J=\int_0^1 t^{\lambda-1} \phi(t)dt+\int_1^\infty  t^{\lambda-1} \phi(t)dt:=j_1+j_2.
\ee
The integral $j_2$ converges classically for any complex $\lambda
\in \C$ and is analytic with respect to $\lambda$.
The integral $j_1$ for $\lambda>0$ converges classically and can be written as
\be\label{e102}
j_1=\int_0^1t^{\lambda-1}(\phi(t)-\phi(0))dt+\phi(0)\frac{t^\lambda|_0^1}{\lambda}=\int_0^1t^{\lambda-1}(\phi(t)-\phi(0))dt+\phi(0)/\lambda.
\ee
The right side of \eqref{e102} admits analytic continuation with respect to $\lambda$ from Re$\lambda>0$ to the region Re$\lambda>-1$. Thus, formulas \eqref{e101} and \eqref{e102}
together define integral \eqref{e100} for Re$\lambda>-1$. 
The singular integral $J$ has a simple pole at $\lambda=0$,
diverges classically for $-1< \lambda<0$, but is defined
in this region by formulas \eqref{e101} and \eqref{e102}
by analytic continuation with respect to $\lambda$.
This procedure can be continued and $J$ can be defined for
an arbitrary large fixed negative $\lambda$, $\lambda\neq 0,-1,-2,...$.

{\bf 1.3.} Let us define the convolution
\be\label{e103}
I(t):=\Phi_\lambda\star b:=\int_0^t\Phi_\lambda(t-s)b(s)ds.
\ee
We assume that $b(t)\in L^1(\R_+)$ and the Laplace transform 
of $b$ is defined for Re$p\ge 0$ by the formula
\be\label{e104}
L(b):=\int_0^\infty e^{-pt}b(t)dt.
\ee
Let us define $L(t^{\lambda-1})$ not using the distribution theory. For $\lambda>0$ one has:
\be\label{e105}
L(t^{\lambda-1})=\int_0^\infty e^{-pt}t^{\lambda-1}dt=
\int_0^\infty e^{-s}s^{\lambda-1}ds p^{-\lambda}=\frac {\Gamma(\lambda)}{p^{\lambda}}.
\ee
It follows from \eqref{e105} and from the definition of $\Phi_{\lambda}=\frac {t^{\lambda-1}}{\Gamma(\lambda)}$ that
\be\label{e106}
L(\Phi_\lambda)=p^{-\lambda}.
\ee
The gamma function $\Gamma(\lambda)$ is analytic in $\lambda\in \C$ except for the simple poles at $\lambda=-n$, $n=0,1,2....$,
with the residue at $\lambda=-n$ equal to $\frac{(-1)^n}{n!}$.
It is known that
\be\label{e107}
\Gamma(z+1)=z\Gamma(z),\quad \Gamma(z)\Gamma(1-z)=\frac \pi {\sin (\pi z)},\quad 2^{2z-1}\Gamma(z)\Gamma(z+\frac 1 2)=\pi^{1/2}\Gamma(2z).
\ee
The function $\frac 1 {\Gamma(\lambda)}$ is an entire function of $\lambda$. 
These properties of $\Gamma(z)$ can be found, for example, in \cite{L}. 

The right side of \eqref{e105} is analytic with respect to $\lambda\in\C$ except for $\lambda=0,-1,-2,....$ and, therefore, defines $L(t^{\lambda-1})$ for all these $\lambda$   by analytic continuation with respect to $\lambda$ without 
using the distribution theory.

Let us define the convolution $I(t)$ using its Laplace transform
\be\label{e108}
L(I(t))=L(\Phi_\lambda)L(b)=\frac {L(b)}{p^{\lambda}}
\ee
and its inverse:
\be\label{e109}
I(t)=L^{-1}\big(L(b)p^{-\lambda}\big),
\ee
where $L^{-1}$ is the inverse of the Laplace transform.
Since the null-space of $L$ is trivial, that is, the zero element, the inverse $L^{-1}$ is well defined on the range of $L$.

For $-1<\lambda<0$, in particular for $\lambda=-\frac 1 4$, formula \eqref{e109}, can be interpreted as a generalized Fourier integral.
The value $\lambda=-\frac 1 4$ is very important in NSP, see monograph \cite{R691}, Chapter 5, and Section 5 below.
We return to this question later, when we discuss the integral equations with hyper-singular integrals. 

{\bf 1.4.} Let us now prove the following result that will be
used later.

{\bf Theorem 1.1.} {\em One has
\be\label{e110}
\Phi_\lambda \star \Phi_\mu=\Phi_{\lambda+\mu}.
\ee
for any $\lambda, \mu\in \C$. If $\lambda+\mu=0$ then
\be\label{e111}
\Phi_0(t)=\delta(t),
\ee
where $\delta(t)$ is the Dirac distribution. 
}

{\em Proof.}  By formulas \eqref{e106} and \eqref{e108} with $b(t)=\Phi_\mu(t)$ one gets
\be\label{e112}
L(\Phi_\lambda\star \Phi_\mu)=\frac 1 {p^{\lambda+\mu}}.
\ee
By formula \eqref{e106} one has 
\be\label{e113}
L^{-1}\Big(\frac 1{p^{\lambda+\mu}}\Big)=\Phi_{\lambda+\mu}.
\ee
This proves formula \eqref{e110}. 

If $\lambda+\mu=0$ then 
\be\label{e113a}
p^{-(\lambda+\mu)}=1, \quad L^{-1}1=\delta(t).
\ee
 This proves formula \eqref{e111}.

Theorem 1.1 is proved. \hfill$\Box$ 

{\bf Remark 1.1.} Let us give an alternative proof of formula \eqref{e110}.  

For Re$\lambda>0$, Re$\mu>0$ one has
\be\label{e114}
\Phi_\lambda\star \Phi_\mu=
\frac 1 {\Gamma(\lambda)\Gamma(\mu)}
\int_0^t(t-s)^{\lambda -1}s^{\mu -1}ds
=\frac{t_+^{\lambda+\mu-1}}{\Gamma(\lambda)\Gamma(\mu)}\int_0^1
(1-u)^{\lambda-1}u^{\mu -1}du=\frac{t_+^{\lambda+\mu -1}}{\Gamma(\lambda+\mu)},
\ee
where the right side of \eqref{e114} is equal to $\Phi_{\lambda+\mu}$ and
 we have used the known formula for the
beta function: 
\be\label{e115}
B(\lambda, \mu):=\int_0^1u^{\lambda -1}(1-u)^{\mu -1}du=
\frac{\Gamma(\lambda)\Gamma(\mu)}
{\Gamma(\lambda+\mu)}.
\ee
Analytic properties of beta function 
follow from these of the gamma function.

{\bf Remark 1.2.} Theorem 1.1 is proved in \cite{GS}, pp.150--151. Our proof differs from the proof in \cite{GS}. It is not clear how the proof in \cite{GS} is related to the definition of regularized hyper-singular integrals used in \cite{GS}.

\section{Preparation for investigation of integral equations\newline with hyper-singular kernels }\label{S:2}

In this section we start an investigation of equations of the following type
 \be\label{e1}
b(t)=b_0(t)+c\int_0^t (t-s)^{\lambda -1}b(s)ds,
\ee 
where  $b_0$ is a smooth functions rapidly decaying with all its derivatives as $t\to \infty$, $b_0(t)=0$ if $t<0$. We are especially interested in the value $\lambda=-\frac 1 4$, because of its importance for the Navier-Stokes theory, \cite{R691}, Chapter 5, \cite{R684}, \cite{R700}.

{\em The integral in \eqref{e1} diverges in the classical sense for $\lambda\le 0$. Our aim is to define this hyper-singular integral.}

 There is a regularization method to define 
singular integrals $J:=\int_{\R}t_+^{\lambda-1}\phi(t)dt$, $\lambda\le 0$, in the distribution theory, see the Introduction, Section 1.2. The integral in \eqref{e1} is a convolution, which is defined in 
\cite{GS}, p.135, as a {\em direct product
of two distributions}.

 This definition
{\em is not suitable} for our purposes because $t_+^{\lambda-1}$ for any $\lambda\le 0$, $\lambda\neq 0,-1,-2,...$ is a distribution 
on the space $\mathcal{K}:=C^\infty_0(\mathbb{R}_+)$ of the test functions, but
{\em it is not a distribution in the
space of the test
functions $K:=C^\infty_0(\R)$ used in \cite{GS}}. 

Indeed, one can find $\phi\in K$ such that
$\lim_{n\to \infty}\phi_n=\phi$ in
 $K$, but $\lim_{n\to \infty}\int_{\R}
t_+^{\lambda-1} \phi_n(t)dt=\infty$
for  $\lambda\le 0$, so that 
$t_+^{\lambda-1}$ is not a linear bounded functional in $K$, i.e., not a distribution.

 For example, the integral
$\int_0^\infty t^{\lambda-1}\phi(t)dt$ is not a bounded linear functional on $K$: take a $\phi$ which is vanishing for $t>1$,  positive near $t=0$ and non-negative on $[0,1]$. Then this integral diverges at such $\phi$ and is not a bounded linear functional on $K$.

On the other hand, one can check that 
$t_+^{\lambda-1}$ for any $\lambda\in R$
is a distribution (a bounded linear functional) in the space 
$\mathcal{K}=C^\infty_0(\R_+)$ with
the convergence $\phi_n\to \phi$ in 
$\mathcal{K}$ defined by the following requirements: 

a) the supports of all $\phi_n$ belong to an interval $[a,b]$,
$0<a\le b<\infty$, 

b) $\phi_n^{(j)}\to 
\phi^{(j)}$ in $C([a,b])$ for all $j=0,1,2,....$. 

Indeed, the functional
$\int_0^\infty t_+^\lambda\phi(t)dt$ is linear and bounded in
$\mathcal{K}$: 
\be\label{e2a}
|\int_0^\infty t_+^\lambda\phi_n(t)dt|\le (a^\lambda+b^\lambda)
\int_a^b |\phi_n(t)|dt.
\ee
A similar estimate holds for all the derivatives of $\phi_n$.

{\em Although $t_+^{-\frac 5 4}$ is a distribution in $\mathcal{K}$, the convolution 
\be\label{e2}
h:=\int_0^t(t-s)^{-\frac 5 4}b(s)ds:=
t_+^{-\frac 5 4}\star b
\ee  
 cannot be defined similarly to the definition in 
the book \cite{GS} because the function\newline
$\int_0^\infty \phi(u+s) b(s)ds$
does not, in general, belong to $\mathcal{K}$ even if 
$\phi\in \mathcal{K}$.}

Let us define the convolution $h$ using the Laplace transform
\eqref{e105}. Laplace transform of distributions is studied in \cite{BP}. There one finds a definition of the Laplace transform of distributions, the Laplace transform of convolutions, tables of the Laplace transforms of distributions, in particular, formula \eqref{e105} and other information.

 One has 
\be\label{e211}
 L(t_+^{-\frac 5 4}\star b)=L(t_+^{-\frac 5 4})L(b).
 \ee
To define $L(t^{\lambda-1})$ for $\lambda\le 0$, note that for 
Re$\lambda>0$ the classical definition \eqref{e105} 
holds. The right side of \eqref{e105}
admits analytic continuation to the
complex plane of $\lambda$, $\lambda\neq 0,-1,-2,....$. This allows one to define integral \eqref{e105}
for any $\lambda\neq 0,-1,-2,...$.
 It is known
that $\Gamma(z+1)=z\Gamma(z)$, so
\be\label{e3a}
\Gamma(-\frac 1 4)=-4\Gamma(3/4):=-c_1, \quad c_1>0. 
\ee
Therefore, {\em we define} 
$h$ by the formula $h=L^{-1}(Lh)$ and defining $L(h)$ as follows:
\be\label{e4}
L(h)=-c_1p^{\frac 1 4}L(b),
\ee
where formula \eqref{e105} with $\lambda=-\frac 1 4$
was used and we assume that $b$ is such that $L(b)$ can be defined. That $L(b)$ is well defined in the Navier-Stokes theory follows from the a priori estimates proved in \cite{R691}, Chapter 5 and in Section 5 below, see Theorem 5.1.

 From \eqref{e4} one gets
\be\label{e4a}
L(b)=-c_1^{-1}p^{-\frac 1 4} L(h).
\ee

\section{Integral equation}\label{S:3}

Consider equation \eqref{e1}. It can be rewritten as
\be\label{e1a}
b(t)=b_0(t)-cc_1\Phi_{\lambda}\star b,
\ee
where
\be\label{e1a1}
 c_1=|\Gamma(-\frac 1 4)|, \quad \lambda=-\frac 1 4.
 \ee
{\bf Theorem 3.1.} {\em Equation \eqref{e1a}-\eqref{e1a1}
has a unique solution in $C(0,T)$ for any $T>0$ if $b_0$
is sufficiently smooth and rapidly decaying as $t$ grows. This solution can be obtained by iterations:}
\be\label{e8} 
b_{n+1}=-(cc_1)^{-1}\Phi_{1/4}\star b_{n} +(cc_1)^{-1}\Phi_{1/4}\star b_0, \quad 
b_{n=0}=
(cc_1)^{-1}\Phi_{1/4}\star b_0, \quad b=\lim_{n\to \infty}b_n.
\ee
{\em Proof.} Applying to  equation\eqref{e1a} the operator $\Phi_{1/4}\star$ and using equation \eqref{e111}
one gets a Volterra-type equation
\be\label{e8a} 
\Phi_{1/4}\star b=\Phi_{1/4}\star b_0-cc_1b,
\ee
or
\be\label{e9}
b=-(cc_1)^{-1}\Phi_{1/4}\star b +(cc_1)^{-1}\Phi_{1/4} \star b_0. 
\ee
The operator $\Phi_\lambda\star$ with 
$\lambda>0$ is a Volterra-type operator. Therefore
equation \eqref{e9} can be solved for $b$ by iterations, see Lemma 3.1 below and \cite{R691}, p.53,
Lemmas 5.10 and 5.11. 

If $b_0\ge 0$ and $cc_1$ is sufficiently large, then
the solution to \eqref{e1a} is non-negative, $b\ge 0$,
see Reamark 3.1 below.
  
Theorem 3.1 is proved.\hfill$\Box$

For convenience of the reader let us 
prove the result about solving equation \eqref{e9} by iterations,  mentioned above.

{\bf Lemma 3.1.} {\em The operator $Af:=\int_0^t(t-s)^pf(s)ds$ in the space 
$X:=C(0,T)$ for any fixed $T\in [0,\infty)$ and $p>-1$ has spectral radius $r(A)$
equal to zero, $r(A)=0$. The equation
$f=Af+g$ is uniquely solvable in $X$. Its solution can be obtained by iterations
\be\label{e9a}
f_{n+1}=Af_n+g, \quad f_0=g; \quad
\lim_{n\to \infty}f_n=f,
\ee
for any $g\in X$ and the convergence 
holds in $X$.}

{\em Proof.}  The spectral radius
of a linear operator $A$ is defined
by the formula $$r(A)=\lim_{n\to \infty}\|A^n\|^{1/n}.$$ By induction one proves that
\be\label{e9b}
|A^nf|\le t^{n(p+1)}\frac{\Gamma^n(p+1)}{\Gamma(n(p+1)+1)}\|f\|_X, \quad n\ge 1.
\ee 
From this formula and the known asymptotic of the gamma function
the conclusion $r(A)=0$ follows.
The convergence result \eqref{e9a}
is analogous to the well known
statement for the assumption $\|A\|<1$. 
A more detailed argument can be found in \cite{R691}, p.53.

Lemma 3.1 is proved. \hfill$\Box$

{\bf Remark 3.1.} {\em If $c>0$ is sufficiently large, then 
the norm of the operator $B:=(cc_1)^{-1}\Phi_{1/4}\star$ in $C(0,T)$
is less than one: $ \|B\|<1$. In this case, $(I-B)^{-1}=\sum_{j=0}^\infty B^j$ is a positive operator.}

Let us now give another approach to solving integral equation \eqref{e1a} with $\lambda=-\frac 1 4$.

{\bf Theorem 3.2.} {\em The solution to equation \eqref{e1a}
with $\lambda=-\frac 1 4$ does exist, is unique, and belongs to $C(\R_+)$ provided that $b_0(t)\in C(\R_+)$ and $ |b_0(t)|+
|b'(t)|\le c(1+t)^{-2}$.}

{\em Proof.} Take the Laplace transform of equation \eqref{e1a}
with $\lambda=-\frac 1 4$, use formula  \eqref{e106} to get
\be\label{20}
L(b)=L(b_0)-cc_1p^{1/4}L(b).
\ee
Thus,
\be\label{21}
L(b)=\frac {L(b_0)}{1+cc_1p^{1/4}}
\ee
Therefore
\be\label{22}
b(t)=L^{-1}\Big(\frac {L(b_0)}{1+cc_1p^{1/4}}\Big).
\ee
Let us check that 
\be\label{23}
max_{t\ge 0}|b(t)|\le c.
\ee
From our assumptions about $b_0(t)$ it follows that $|L(b_0)|
\le c(1+|p|)^{-1}$, Re$p\ge 0$. Let $p=iw$.
Since $b(t)=(2\pi)^{-1}\int_{-\infty}^{\infty}e^{iwt}L(b)dw$, one gets
\be\label{24}
|b(t)|\le\frac c {2\pi}\int_{-\infty}^{\infty}(1+|w|)^{-1}
|1+cc_1(iw)^{1/4}|^{-1}<c_2,
\ee
where $c_2>0$ is some constant.
Here we have used the inequality
\be\label{24a}
\inf_{w\in \R}|1+cc_1(iw)^{1/4}|^{-1}\le c.
\ee
Recall that by $c>0$ various constants are denoted.

Let us check \eqref{24a} for $w\ge 0$. For $w<0$ the argument is similar. One has $(iw)^{1/4}=e^{i\pi/8}w^{1/4}$, $$J:=\frac 1 {|1+C
\cos(\pi/8)+iC\sin(\pi/8)|},$$
where $C:=cc_1w^{1/4}>0$. Therefore,
$$J^{-2}=[1+C\cos(\pi/8)]^2+C^2\sin^2(\pi/8)=1+C^2+2C\cos(\pi/8)>c>0.$$
Consequently, inequality \eqref{24a}
is checked.

Theorem 3.2 is proved. \hfill$\Box$

{\bf Remark 3.2.} {\em It follows from formula \eqref{e9}
that $b(0)=0$ because $\lim_{t\to 0}\Phi_{\frac 1 4}\star b_0=0$
and $\lim_{t\to 0}\Phi_{\frac 1 4}\star b=0$ holds if $b$ is a locally integrable bounded on $\R_+$ function $b=b(t)$.}

\section{Integral inequality}\label{S:4}
Consider the following inequality
\be\label{e7}
q(t)\le b_0(t)+ct_+^{\lambda-1}\star q=b_0(t)-cc_1\Phi_{-\frac 1 4}q,
\ee
where $c_1=-\Gamma(-\frac 1 4)$ for $\lambda=-\frac 14$

Let $f=f(t)\in L^1(\R_+)$ be some function. If
$q\le f$ then $\Phi_{1/4}\star q\le \Phi_{1/4}\star f.$ 
Therefore,  inequality \eqref{e7} with $\lambda=-\frac 1 4$, after applying to both sides the operator $\Phi_{1/4}\star$, implies 
\be\label{e9c}
q\le -(cc_1)^{-1}\Phi_{1/4}\star q +(cc_1)^{-1}\Phi_{1/4}\star b_0.
\ee
Inequality \eqref{e9c} for sufficiently large $c>0$ can be solved by iterations with the initial term
$(cc_1)^{-1}\Phi_{1/4}\star b_0$, see Remark 3.1. This
yields
\be\label{e11}
q(t)\le b(t),
\ee 
where $b$ solves the integral equation \eqref{e1a}.
This follows from Theorem 4.1.

{\bf Theorem 4.1.} {\em Assume that $b$ solves \eqref{e1a}, $c>0$ is sufficiently large, $b_0(t)$ satisfies conditions stated in Theorem 3.2 and $q\ge 0$
solves inequality \eqref{e7}. Then inequality \eqref{e11} holds.}

{\em Proof.} Denote $z:=b-q$, where 
\be\label{e9cd}
b=-(cc_1)^{-1}\Phi_{1/4}\star q +(cc_1)^{-1}\Phi_{1/4}\star b_0.
\ee
Then
\be\label{e11a}
0\le z+(cc_1)^{-1}\Phi_{\frac 1 4}\star z.
\ee
Solving this inequality by iterations and using Remark 3.1 one obtains \eqref{e11}. If $c>0$ is arbitrary, then this argument yields inequality \eqref{e11} for sufficiently small $t>0$ because the norm of the operator $(cc_1)^{-1}\Phi_{1/4}\star$ tends to zero when $t\to 0$.

Theorem 4.1 is proved. \hfill$\Box$

Papers \cite{R677}, \cite{R698}, \cite{R704} also deal with hyper-singular integrals.

\section{Application to the Navier-Stokes problem}\label{S:5}

In this Section we apply the results of Sections 1--4
to the Navier-Stokes problem.  Especially the results of
Sections 3 and 4 will be used.

The Navier-Stokes problem 
(NSP) in $\R^3$ is discussed in many books and papers ( see \cite{R691}, Chapter 5, and references therein).
 The uniqueness of a solution in $\R^3$ was proved in \cite{La}, \cite{R691} and in \cite{R700} in different norms. 
 The existence of the solution to the NSP is discussed
 in \cite{R691}.  
  
The goal of this Section is to prove that the statement of the NSP is contradictory. Therefore, the NSP is not a
physically correct statement of the problem of fluid mechanics.
We prove that the solution to the NSP does not exist,
in general. Therefore, in this Section a negative solution to one of the millennium problems is given.

{\em  What is a physically correct statement of problems
of fluid mechanics remains an open problem.}

We prove {\em the paradox in the NSP}. This paradox  can be described as follows:

{\em One can have initial velocity $v(x,0)>0$ in the NSP and, nevertheless, the
solution $v(x,t)$ to this NSP must have the zero initial velocity: $v(x,0)=0$. }

This paradox proves that the statement of the NSP is contradictory, that the NSP is not a physically correct statement of the fluid mechanics problem and the solution to the NSP does not exist, in general.

 The NSP in $\R^3$ consists of solving the equations
\be\label{e501} v'+(v, \nabla)v=-\nabla p +\nu\Delta v +f, \quad x\in \R^3,\,\, t\ge 0,\quad \nabla \cdot v=0,\quad v(x,0)=v_0(x),
\ee
see, for example, books \cite{La} and \cite{R691}, Chapter 5.

Vector-functions velocity $v=v(x,t)$ and the exterior force $f=f(x,t)$ and the scalar function $p=p(x,t)$,
the pressure, are assumed to decay as $|x|\to \infty$ and
$t\in \R_+:=[0, \infty)$.

 The derivative with respect to time is denoted $v':=v_t$,
$\nu=const>0$ is the viscosity coefficient,  the velocity $v=v(x,t)$ and the pressure $p=p(x,t)$ are unknown,  $v_0=v(x,0)$ and $f(x,t)$ are known. It is assumed that  $ \nabla \cdot v_0=0$.
Equations \eqref{e501}  describe viscous incompressible fluid with density $\rho=1$.

Let us assume for simplicity that $f=0$. This do not change our arguments and our logic.

 The solution to NSP \eqref{e501} solves the  integral equation:
\be\label{e502} v(x,t)=F- \int_0^tds \int_{\R^3} G(x-y,t-s)(v,\nabla)v dy, 
\ee
where $(v,\nabla)v=v_j v_{p,j}$, over the repeated indices summation is assumed and 
$v_{p,j}:=\frac {\partial v_p}{\partial x_j}$.

Equation \eqref{e502} implies an integral inequality
of the type studied in Sections 3 and 4 (see also \cite{R691}, Chapter 5).

 Formula for the tensor $G=G(x,t)=G_{pm}(x,t)$
is derived in \cite{R691}, p.41: 
\be\label{e5021}
G(x,t)=(2\pi)^{-3}\int_{\R^3} e^{i\xi \cdot x}\Big( \delta_{pm}-\frac {\xi_p \xi_m}{\xi^2}\Big)e^{-\nu \xi^2 t}d\xi.
\ee
The term $F=F(x,t)$,  in our case when $f=0$, depends only on the data $v_0$ (see  formula (5.42) in \cite{R691}):
\be\label{e503}
F(x,t):=\int_{\R^3}g(x-y,t)v_0(y)dy,
\ee
where
\be\label{e503'}
g(x,t)=\frac{e^{-\frac{|x|^2}{4\nu t}}}{4\nu \pi t}, \quad t>0;
\quad g(x,t)=0, \quad t\le 0.
\ee
We assume throughout that 
\be\label{e503"}
v_0=v(x,0)>0
\ee 
is such that 
$F$ is bounded in all the norms we use.

Let us use the Fourier transform: 
\be\label{e5044}
\tilde{v}:=\tilde{v}(\xi,t):=(2\pi)^{-3}\int_{\R^3}v(x,t)e^{-i\xi \cdot x}dx.
\ee
 Fourier transform  equation \eqref{e502} and get the integral equation:
\be\label{e504}
\tilde{v}(\xi,t)=\tilde{F}(\xi,t)-\int_0^tds \tilde{G}(\xi,t-s) \tilde{v}\bigstar (i\xi \tilde{v}), 
\ee
where  $\bigstar$ denotes the convolution in $\R^3$. 

 For brevity we omitted the tensorial  indices:
instead of $\tilde{G}_{mp}\tilde{v}_j\bigstar (i\xi_j)\tilde{v}_p$, where one sums up over the repeated indices, we wrote
$ \tilde{G} \tilde{v}\bigstar (i\xi \tilde{v})$.

From formula (5.9) in \cite{R691}, see formula \eqref{e5021} one gets: 
\be\label{e5050}
\tilde{G}(\xi,t)=(2\pi)^{-3}\Big( \delta_{pm}-\frac {\xi_p \xi_m}{\xi^2}\Big)e^{-\nu \xi^2 t}.
\ee
One has $|\delta_{pm}-\frac {\xi_p \xi_m}{\xi^2}|\le c$.
Therefore,
\be\label{e505}
|\tilde{G}(\xi,t-s)|\le ce^{-\nu (t-s) \xi^2}.
\ee 
 We denote by $c>0$ {\em various constants} independent of $t$ and $\xi$, by $\|\tilde{v}\|$ the norm in $L^2(\R^3)$ and by $(v,w)$ the inner product in $L^2(\R^3)$.
 
 Let us introduce the norm
\be\label{e505a}
\|v\|_1:=\|v\|+\|\nabla v\|. 
\ee
  One has
\be\label{e505b} (2\pi)^{3/2}\|\tilde{v}\|= \|v\|, \quad (2\pi)^3\||\xi|\tilde{v}\|^2=\|\nabla v\|^2,
 \ee
by the Parceval equality. 
  
 {\bf Assumption A.}  {\em Assume that $F(x,t)$ is a smooth function rapidly decaying together with all its derivatives. In particular, 
$$\sup_{t\ge 0}\Big((1+t^m)\|F(x,t)\|_1\Big)+\sup_{t\ge 0, \xi\in \R^3}\left((1+t^m+|\xi|^m)|\tilde{F}(\xi,t)|\right)<c, \quad  m=1,2,3.$$ }

Assumption A holds throughout Section 5 and is not repeated. It is known  that
\be\label{e507a}
\sup_{t\ge 0}\left(\|v\|+\int_0^t \|\nabla v\|^2ds\right)<c, \quad \sup_{t\ge 0}\int_0^t \|\tilde{v}|\xi|\|^2ds<c,
\ee
 \be\label{e507} 
\sup_{t\ge 0}(|\xi||\tilde{v}(\xi,t)|)<\infty,
\quad |\tilde{v(\xi,t)}|\le c(1+t^{1/2}), \quad \sup_{t\ge 0}\|\nabla v\|< \infty,
\ee 
see \cite{R691}, p.52.

{\bf Theorem 5.1.} {\em Inequalities \eqref{e507a}--\eqref{e507} hold.  }

{\bf Theorem 5.2.} {\em The NSP \eqref{e501} does not have a solution, in general.}
 
{\bf Proof of Theorem 5.1.} Proof of Theorem 5.1 can be found in \cite{R691}.  Because of the importance
 of the third inequality \eqref{e507} and of its novelty, we give its proof
 in detail.

 Let $|\tilde{v}(\xi,t)|:=u$, $|\tilde{F}|:=\mu(\xi,t):=\mu$. From equation \eqref{e504} one gets:
\be\label{e514a}
u\le \mu+c\int_0^t e^{-\nu (t-s)\xi^2}\|u\|\||\xi|u\|ds\le \mu+c\int_0^t e^{-\nu (t-s)\xi^2}b(s)ds, \quad b(s):=\||\xi|u\|,
\ee
where the Parceval formula
\be\label{e51410}
(2\pi)^{3/2} \|\tilde{v}\|=\|v\|<c 
\ee
was used. 

By direct calculation one derives the following inequality:
\be\label{e51411}
\| e^{-\nu (t-s)\xi^2}|\xi|\|\le c(t-s)^{-\frac 5 4}.
\ee
It follows from this inequality and from \eqref{e514a} by multiplying by $|\xi|$
and taking the norm $\|\cdot\|$ of the resulting inequality that the following integral inequality holds:
\be\label{e514}
b(t)\le b_0(t)+c\int_0^t (t-s)^{-\frac 5 4}b(s)ds,
 \ee 
where
\be\label{e51412}
b_0(t):=\||\xi|\mu(\xi,t)\|,
 \quad b(s):=\||\xi|u\|.
 \ee
The function $b_0(t)$ is smooth and rapidly decaying due to Assumption A.

Let $\beta$ solve the following equation:
\be\label{e514b}
\beta(t)=b_0(t)+c\int_0^t (t-s)^{-\frac 5 4}\beta(s)ds.
\ee
Equation \eqref{e514b} can be written as
\be\label{e514bb}
\beta(t)=b_0(t)-cc_1\Phi_{-\frac 14}\star \beta,
\ee
where $\star$ denotes the convolution of two functions on $\R_+$ and $c_1=|\Gamma(-\frac 1 4)|$. The convolution
on $\R_+$ was defined in the Introduction. In Section 4 the relation between the solutions to integral equation  \eqref{e514bb} and integral inequality \eqref{e514} was
studied and the inequality \newline
$b(t)\le \beta(t)$ was proved.

Taking the Laplace transform of equation \eqref{e514b} and using 
equation \eqref{e105}
 with $\lambda=-\frac 1 4$, we get 
\be\label{e514cc}
L(\Phi_{-\frac 14}\star \beta)=L(\Phi_{-\frac 14})L(\beta)=p^{1/4} L(\beta),
\ee
so 
\be\label{e514ccc}
L(\beta)=L(b_0)-cc_1p^{1/4}L(\beta).
\ee
Therefore,
\be\label{e514c}
L(\beta)=\frac {L(b_0)}{1+cc_1p^{1/4}}, \quad 0\le b(t)\le \beta(t).
\ee
It follows from \eqref{e514c} that
\be\label{e514d}
b(t)\le \beta(t)=\frac 1{2\pi}\int_{-\infty}^{\infty}
e^{i\tau t}\frac {L(b_0)}{1+cc_1(i\tau)^{1/4}}d\tau\le \frac 1{\pi}\int_0^\infty
\frac {|L(b_0)|}{|1+cc_1e^{i\pi/8}\tau^{1/4}|}d\tau\le c,
\ee
where the argument $p$ of the function $L(b_0)$ is equal to $i\tau$,
$p=i\tau$, and we have used the decay
$O(|\tau|^{-1})$ of $|L(b_0)|$ as a function of 
 $p=i\tau$ as $|\tau|\to \infty$.
 
  This decay
follows from Assumption A and implies that
the integrand  in \eqref{e514d} belongs to $L^1(\R)$ 
because of the following inequality, proved at the end of Section 3:
\be\label{e6111} 
\inf_{\tau\in [0,\infty)}|1+cc_1e^{i\pi/8}\tau^{1/4}|>0,
\quad cc_1>0.
\ee 
Thus, the third estimate \eqref{e507} of Theorem 5.1  is proved.\hfill$\Box$

 {\em Proof of theorem 5.2.} If $v_0(x)=v(x,0)\not\equiv 0$ and $\nabla \cdot v_0(x)=0$ then 
 $b_0(0)>0$. Apply to equation \eqref{e514b} the operator $\Phi_{1/4}\star$ and use Theorem 1.1. This yields
\be\label{e61} 
\Phi_{1/4}\star \beta=\Phi_{1/4}\star b_0-cc_1\beta(t),
\ee
where formula \eqref{e105} was used, $c_1=-\Gamma(-\frac 1 4)>0$ and $\Phi_{\frac 1 4}\star\Phi_{-\frac 1 4}=\delta$, where $\delta$ is the delta-function, see formulas  \eqref{e110}--\eqref{e111}.  We assume that $b_0(t)$ satisfies Assumption A, so it is smooth and rapidly decaying. Then equation \eqref{e61} can be solved by iterations by Theorem 3.1
and the solution $\beta$ is also smooth. Therefore,
\be\label{e71} 
 \lim_{t\to 0}\Phi_{1/4}\star \beta=0, \quad \lim_{t\to 0}\Phi_{1/4}\star b_0=0.
 \ee
 Consequently, it follows from \eqref{e61} that 
  $\beta(0)=0$. 
  Since $0\le b(t)\le \beta(t)$, one concludes that
\be\label{e72} 
  b(0)=0.
  \ee

This result proves that the NSP problem  does not have a solution, in general. Indeed, starting with an initial data  which is positive we prove that the corresponding solution to the NSP must have the initial data equal to zero. This is the {\em NSP paradox},
see \cite{R708}. 

Of course, if the data $v_0(x)=v(x,0)=0$
then the solution exists for all $t\ge 0$ and is equal to zero by the uniqueness theorem, see, for example, \cite{R691}, \cite{R700}.

Other paradoxes of the theory of fluid mechanics are mentioned in \cite{La}.

Theorem 5.2 is proved. \hfill$\Box$

\end{document}